\documentclass{elsart3-1}

 \usepackage{graphics, ulem}

 \usepackage{epsfig}
\usepackage{color}
\usepackage{amssymb}

\usepackage[english,francais]{babel}

\newtheorem{e-proposition}[theorem]{Proposition}

\newtheorem{e-definition}[theorem]{Definition\rm}

\newtheorem{theoreme}{Th\'eor\`eme}[section]
\newtheorem{lemme}[theoreme]{Lemme}

\newtheorem{remarque}{\it Remarque}

\renewcommand{\theequation}{\arabic{equation}}
\def\R {\mathbf{R}}

\def\T {\mathbf{T}}

\def\Z {\mathbf{Z}}

\def\cD {\mathcal{D}}

\def\cZ {\mathcal{Z}}

\def\eps {{\varepsilon}}

\def\indc {{\bf 1}}

\def\d {{\partial}}

\def\og{\leavevmode\raise.3ex\hbox{$\scriptscriptstyle\langle\!\langle$~}}
\def\fg{\leavevmode\raise.3ex\hbox{~$\!\scriptscriptstyle\,\rangle\!\rangle$}}

\journal{the Acad\'emie des sciences}
\begin{document}
% Vous pouvez mettre dans la prochain ligne la rubrique choisie
% (si vous la connaissez) - meme deux, format : Rubrique1/Rubrique2
\centerline{}
\begin{frontmatter}

\selectlanguage{francais}
\title{Limite de diffusion lin\'eaire \\pour un syst\`eme d\'eterministe de sph\`eres dures}

% utiliser les \`etiquettes pour indiquer l'adresse de chaque auteur,
%     s'il y a plusieurs adresses

% \author[label1,label2]{}
% \address[label1]{}
% \address[label2]{}

\author[thierry]{T. Bodineau},
\ead{Thierry.Bodineau@ens.fr}
\author[isabelle]{I. Gallagher}
\ead{gallagher@math.univ-paris-diderot.fr}
\author[laure]{L. Saint-Raymond}
\ead{Laure.Saint-Raymond@ens.fr}

\address[thierry]{CNRS \& \'Ecole Normale Sup\'erieure, D\'epartement de Math\'ematiques et Applications}
\address[isabelle]{Universit\'e Paris-Diderot, Institut de Math\'ematiques de Jussieu}
\address[laure]{Universit\'e Pierre et Marie Curie \& \'Ecole Normale Sup\'erieure, D\'epartement de Math\'ematiques et Applications}

\medskip
\selectlanguage{francais}
\begin{center}
{\small Re\c{c}u le *****~; accept\'e apr\`es r\'evision le +++++\\
Pr\'esent\'e par Jean-Pierre Kahane}
\end{center}

\begin{abstract}
\selectlanguage{francais}
Cette note montre comment on peut obtenir le mouvement brownien   comme limite hydrodynamique d'un  syst\`eme d\'eterministe de sph\`eres dures quand le nombre de particules $N$ tend vers l'infini et que leur diam\`etre~$\eps$ tend vers~0, dans la limite de relaxation rapide $N\eps^{d-1}\to \infty $ (avec un choix d'\'echelles de temps et d'espace convenable).
Comme sugg\'er\'e par Hilbert dans son sixi\`eme probl\`eme, on utilise  la th\'eorie cin\'etique de Boltzmann comme niveau de description  interm\'ediaire.
La preuve suit les id\'ees fondamentales de Lanford sur la propagation du chaos. La principale nouveaut\'e consiste \`a obtenir des estimations sur les arbres de collision pathologiques par une \'etude fine du processus de branchement.

%{\it Pour citer cet article~: A. Nom1, A. Nom2, C. R.
%Acad. Sci. Paris, Ser. I 340 (2005).}
\vskip 0.5\baselineskip

\selectlanguage{english}
\noindent{\bf Abstract}
\vskip 0.5\baselineskip
\noindent
{\bf Linear diffusive limit of deterministic systems of hard spheres. }
We provide a rigorous derivation of the brownian motion    as the hydrodynamic limit of a deterministic system of hard-spheres as the number of particles~$N$ goes to infinity and their diameter $\eps$ simultaneously goes to $0,$ in the fast relaxation limit $N\eps^{d-1}\to \infty $ (with a suitable scaling of the observation time and length).
As suggested by Hilbert in his sixth problem, we use   Boltzmann's kinetic theory as an intermediate level of description for  the gas close to global equilibrium.
Our proof relies on the fundamental ideas of Lanford. 
The main novelty  is the detailed study of the branching process, leading to   explicit estimates on pathological collision trees.

%{\it To cite this article: A. Nom1, A. Nom2, C. R.
%Acad. Sci. Paris, Ser. I 340 (2005).}
\end{abstract}
\end{frontmatter}

\selectlanguage{francais}
% texte principale
\section{Introduction}
\label{}

Le sixi\`eme probl\`eme pr\'esent\'e par Hilbert au Congr\`es International des Math\'ematiciens en 1900 concerne l'axiomatisation de la physique, et plus pr\'ecis\'ement la description math\'ematique de la consistance  entre les mod\`eles atomiques et les mod\`eles continus de la dynamique des gaz.
De fa\c con un peu plus restrictive (puisque seuls les gaz parfaits peuvent \^etre appr\'ehend\'es dans cette approche), Hilbert sugg\`ere d'utiliser l'\'equation cin\'etique de Boltzmann comme \'etape interm\'ediaire pour comprendre l'apparition de l'irr\'eversibilit\'e  et des m\'ecanismes de diffusion \cite{hilbert}:

`` {\it Quant aux principes de la M\'ecanique, nous poss\'edons d\'ej\`a au point de vue physique des recherches d'une haute port\'ee; je citerai, par exemple, les \'ecrits de MM. Mach, Hertz, Boltzmann et Volkmann. Il serait aussi tr\`es d\'esirable qu'un examen approfondi des principes de la M\'ecanique f\^ut alors tent\'e par les math\'ematiciens. Ainsi le Livre de M. Boltzmann sur les Principes de la M\'ecanique nous incite \`a \'etablir et \`a discuter au point de vue math\'ematique d'une mani\`ere compl\`ete et rigoureuse les m\'ethodes bas\'ees sur l'id\'ee de passage \`a la limite, et qui de la conception atomique nous conduisent aux lois du mouvement des continua. Inversement on pourrait, au moyen de m\'ethodes bas\'ees sur l'id\'ee de passage \`a la limite, chercher \`a d\'eduire les lois du mouvement des corps rigides d'un syst\`eme d'axiomes reposant sur la notion d'\'etats d'une mati\`ere remplissant tout l'espace d'une mani\`ere continue, variant d'une mani\`ere continue et que l'on devra d\'efinir param\'etriquement.

Quoi qu'il en soit, c'est la question de l'\'equivalence des divers syst\`emes d'axiomes qui pr\'esentera toujours l'int\'er\^et le plus grand quant aux principes.}"

\medskip
Une litt\'erature importante est consacr\'ee \`a ces probl\`emes d'analyse asymptotique, mais ils restent \`a ce jour encore tr\`es ouverts. 

Les travaux fondamentaux de DiPerna et Lions sur les solutions renormalis\'ees de l'\'equation de Boltzmann   \cite{diperna-lions} ont permis  d'obtenir une \'etude compl\`ete de certaines limites hydrodynamiques de l'\'equation de Boltzmann, en particulier dans le  r\'egime incompressible visqueux conduisant aux \'equations de Navier-Stokes (\cite{BGL2,GSR}).

A ce stade, le principal obstacle semble donc \^etre li\'e \`a l'autre \'etape, c'est-\`a-dire \`a l'obtention  de l'\'equation de Boltzmann \`a partir de syst\`emes de particules : le r\'esultat le plus g\'en\'eral concernant cette limite de faible densit\'e, qui est d\^u \`a Lanford dans le cas des sph\`eres dures \cite{lanford} (voir aussi  \cite{CIP,uchiyama,GSRT}, pour une preuve compl\`ete)  n'est en effet valable que pour des temps tr\`es courts, ce qui ne permet pas d'observer une quelconque relaxation.

\subsection{L'\'equation de Boltzmann, un mod\`ele interm\'ediaire}

L'\'etat d'un gaz monoatomique (constitu\'e de particules indiscernables) peut \^etre caract\'eris\'e par sa fonction de distribution $f$,  o\`u $f(t,x,v)$ est la densit\'e associ\'ee \`a la probabilit\'e de trouver une particule de position 
$x$ et de vitesse $v$ \`a l'instant $t$.
Cette fonction ne se mesure pas, mais permet de calculer les observables (telles que la temp\'erature $\Theta$ ou la vitesse d'\'ecoulement $U$) en prenant des moyennes~:
\begin{eqnarray*}
R(t,x) &:= &\int f(t,x,v) dv \, ,\quad 
 U(t,x):=\frac1{R(t,x)} \int v f(t,x,v) dv\, ,\\
  \Theta(t,x)& := &\frac1{3R(t,x)} \int f(t,x,v) |v-U(t,x)|^2 dv \,.
\end{eqnarray*}

Dans le vide et en l'absence d'interactions, les particules ont un mouvement rectiligne et uniforme. 
La distribution $f$ satisfait alors l'\'equation de transport libre
$$\d_t f+v \cdot \nabla_x f = 0\,. $$
Ici on va supposer que les particules - qui sont de taille n\'egligeable -  interagissent  par contact, et que ces collisions sont \'elastiques (i.e. pr\'eservent la quantit\'e de mouvement et l'\'energie). Une analyse tr\`es similaire mais techniquement plus difficile permettrait en fait de consid\'erer des interactions hamiltoniennes d\'ecrites par un potentiel \`a support compact (\cite{king,GSRT}).
Sous l'effet combin\'e du transport et des collisions, la distribution $f $  des particules \'evolue alors selon l'\'equation int\'egro-differentielle suivante, dite {\bf \'equation de Boltzmann non lin\'eaire}
\begin{equation}
\label{NL-boltz}
\begin{array}l
 \displaystyle \d_t f +v\cdot \nabla_x f=Q(f,f) \, ,\\
 \displaystyle Q(f,f)(v):=\int \Big( f(v')f(v'_1)-f(v)f(v_1)\Big)   \,((v-v_1)\cdot \omega)_+ \, dv_1 d\omega\,.
  \end{array}
\end{equation}
 Comme les collisions sont \'elastiques, les vitesses pr\'e-collisionnelles $(v',v'_1)$ v\'erifient 
  $$
  v'+v'_1 = v+v_1,\qquad  |v'|^2+|v'_1|^2 = |v|^2+|v_1|^2,
    $$
et peuvent \^etre param\'etris\'ees par l'angle de d\'eflection  $\omega\in {\mathbb S}^{d-1}$ 
$$v'_1 := v_1 +\Big( (v-v_1)\cdot \omega\Big)  \,\omega, \quad v':=v-\Big( (v-v_1)\cdot \omega\Big) \, \omega\,.$$
\begin{figure}[h] %  figure placement: here, top, bottom, or page
   \centering
  \includegraphics[width=2in]{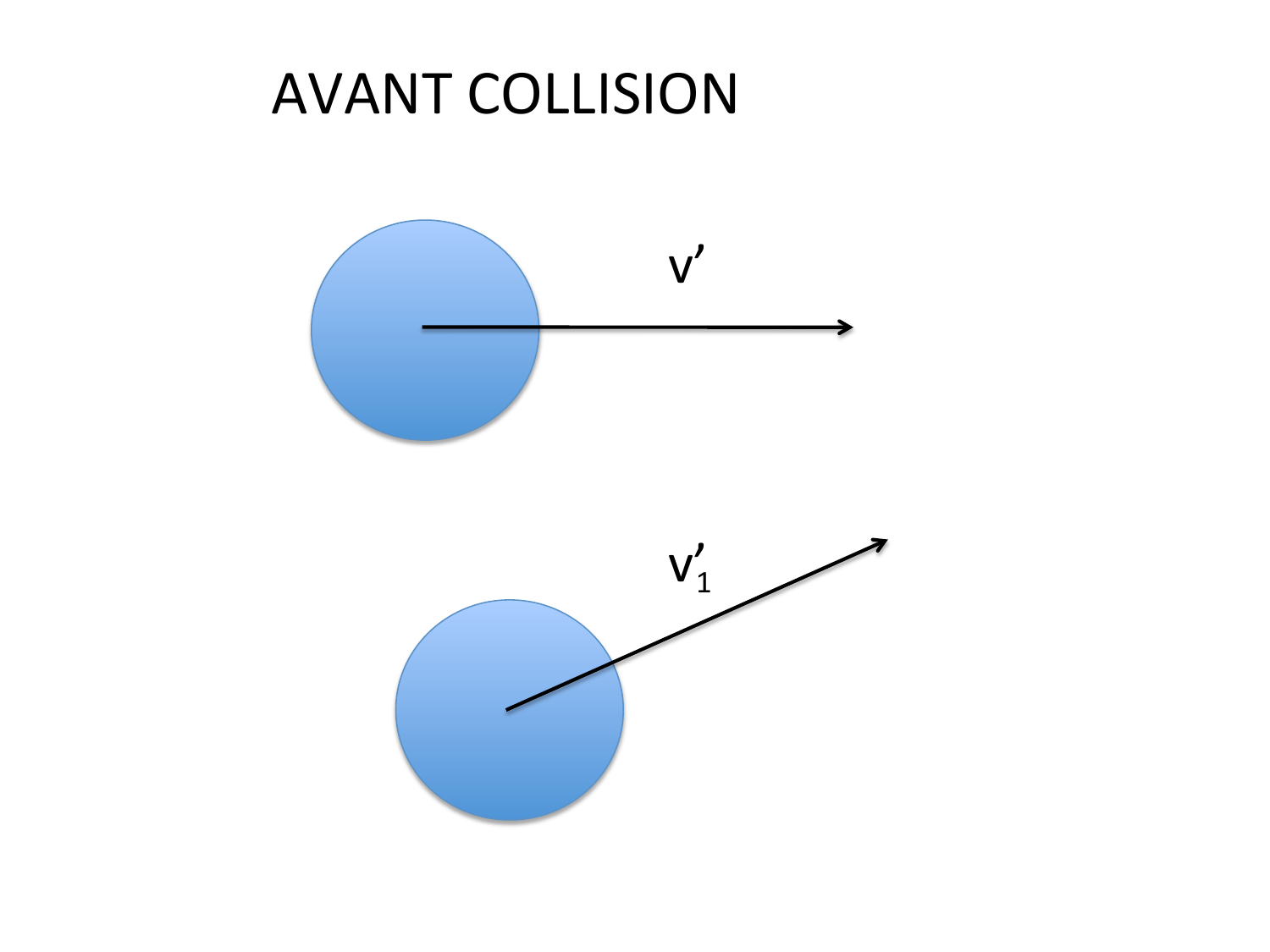} \includegraphics[width=2in]{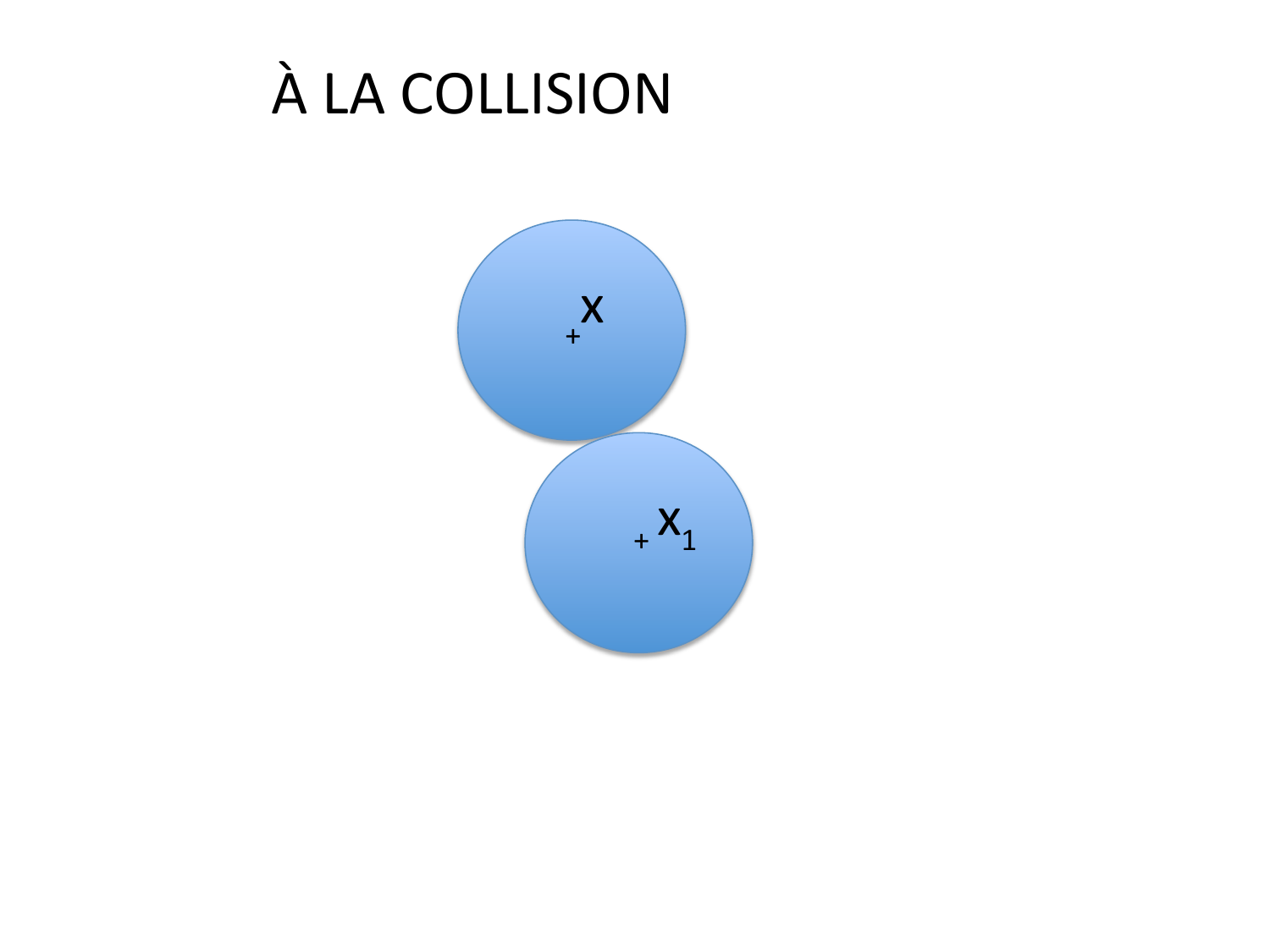} \includegraphics[width=2in]{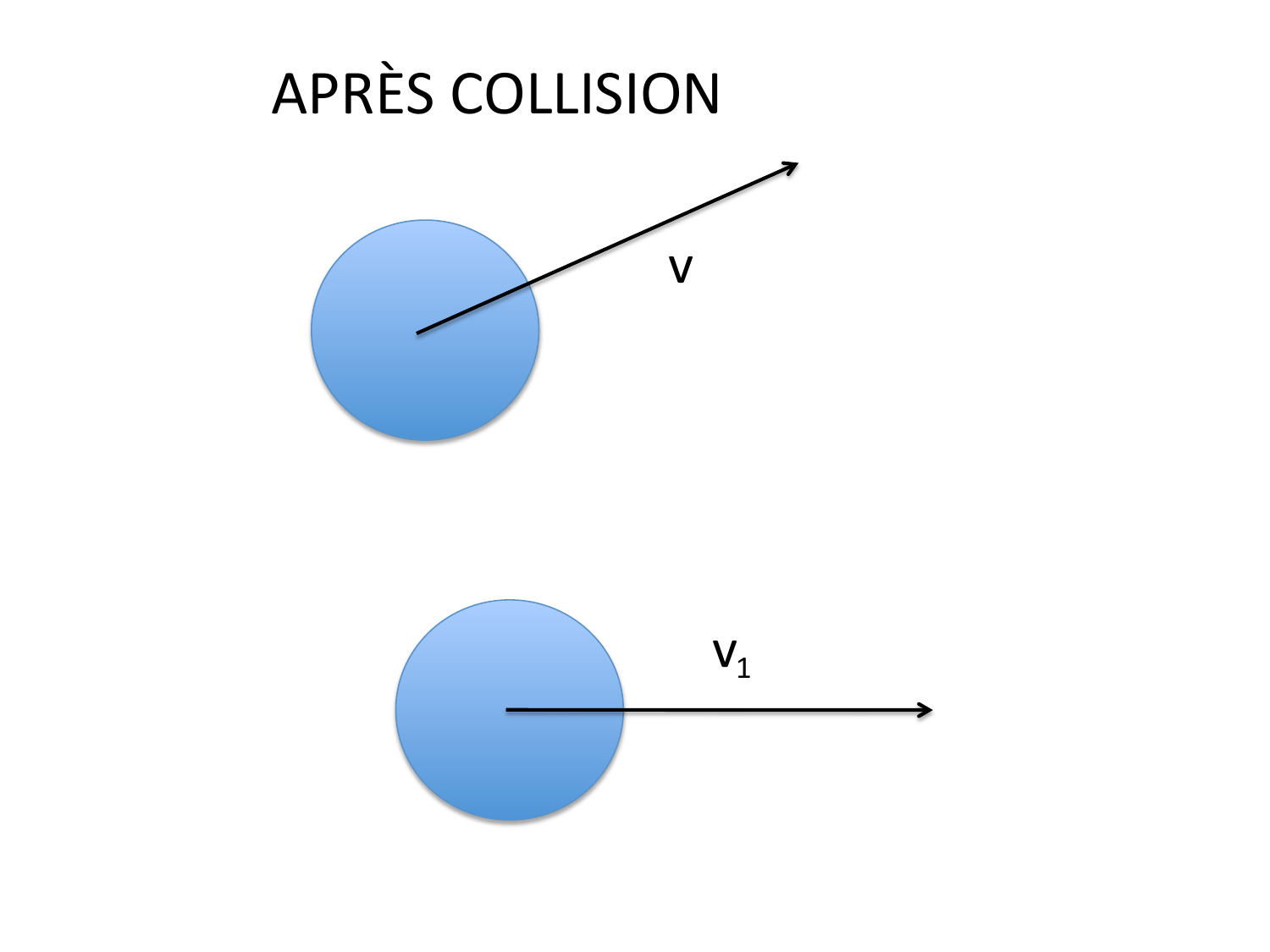} 
   \caption{Collision de sph\`eres dures}
 \end{figure}

Ë cause de l'indiscernabilit\'e des particules, $v$ et $v_1$ jouent des r\^oles sym\'etriques dans l'int\'egrale de collision. La r\'eversibilit\'e du m\'ecanisme \'el\'ementaire de collision  implique de plus que le changement de variables $(v',v'_1)\to (v,v_1)$ est de jacobien 1, de sorte que
$$\int Q(f,f) \varphi(v) dv =\frac14 \int [f(v')f(v'_1)-f(v)f(v_1)] \big( \varphi(v)+\varphi(v_1)-\varphi(v')-\varphi(v'_1)\big)  \, \big((v-v_1)\cdot \omega\big)_+ \, dv dv_1 d\omega \,.$$
En particulier, en choisissant $\varphi(v) = 1$, puis $\varphi(v) = v$ et $\varphi(v) =|v|^2$, on obtient formellement les conservations de la masse, de l'impulsion et de l'\'energie
$$
\begin{array}l
\displaystyle \d_t R +\nabla_x \cdot (R U) =0\,,\\
\displaystyle \d_t (R U) +\nabla_x \cdot \int f v\otimes v dv =0\,,\\
\displaystyle \d_t (R U^2 +3R\Theta ) +\nabla_x \cdot \int f |v|^2 vdv =0\,.
\end{array}
$$

\subsection{Entropie, irr\'eversibilit\'e et relaxation vers l'\'equilibre}

En prenant $\varphi =\log f$ dans l'identit\'e pr\'ec\'edente, on obtient aussi
$D(f) :=-\int Q(f,f) \log f(v) \geq 0,$
d'o\`u l'on d\'eduit l'in\'egalit\'e d'entropie, aussi appel\'ee th\'eor\`eme H de Boltzmann,
$$
\displaystyle \int f\log f(t,x,v) dxdv+\int_0^t \int D(f) (s,x) dsdx \leq \int f_0\log f_0(x,v) dxdv\,.
$$
Cela signifie que l'\'equation de Boltzmann d\'ecrit une {\bf \'evolution irr\'eversible}.

Les maxima de l'entropie \`a masse, impulsion et \'energie fix\'ees 
sont donn\'es formellement par le th\'eor\`eme  d'optimisation sous contrainte de Lagrange 
$$\log f (v) =\alpha +\beta \cdot v +\gamma |v|^2\,.$$
En d'autres termes, les distributions d'\'equilibre sont  des Gaussiennes, ce qui est conforme \`a la pr\'ediction de Maxwell.
On s'attend \`a ce que  l'\'equation de Boltzmann pr\'edise une {\bf relaxation vers ces \'equilibres}.

Ces deux caract\'eristiques de la dynamique de Boltzmann montrent que la description statistique n'est pas tout-\`a-fait \'equivalente \`a la description microscopique. En effet,
\begin{itemize}
\item l'entropie  cro\^\i t alors que les \'equations de la m\'ecanique sont r\'eversibles (paradoxe de Lodschmidt);
\item la relaxation est incompatible avec le th\'eor\`eme de Poincar\'e qui pr\'edit une r\'ecurrence pour le syst\`eme de $N$ sph\`eres dures (paradoxe de Zermelo)\,.
\end{itemize}
Obtenir  l'\'equation de Boltzmann \`a partir de la m\'ecanique de Newton n\'ecessite donc de comprendre l'origine de l'irr\'eversibilit\'e, irr\'eversibilit\'e qui est aussi cruciale pour les limites de diffusion.

\section{R\'esultats de convergence}

Les travaux de Grad montrent que la  mesure empirique associ\'ee au syst\`eme de $N$ sph\`eres satisfait une \'equation collisionnelle similaire \`a l'\'equation de Boltzmann, c'est-\`a-dire o\`u le transport et les collisions sont deux ph\'enom\`enes qui s'\'equilibrent,  si 
 le nombre de particules $N$ tend vers l'infini,
et leur diam\`etre $\eps$ tend simultan\'ement vers 0, avec le choix d'\'echelles de Boltzmann-Grad $N\eps^{d-1} \sim1\,.$
Le point \`a noter est que l'op\'erateur de collision porte  sur la probabilit\'e jointe de trouver deux particules collisionnelles, qui en g\'en\'eral ne s'exprime pas en fonction de la densit\'e de probabilit\'e \`a une particule.

La validit\'e de l'\'equation de Boltzmann est donc li\'ee \`a une propri\'et\'e d'ind\'ependance des particules collisionnelles, que l'on appelle souvent ``chaos".
Le th\'eor\`eme de Lanford \cite{lanford} (voir aussi \cite{GSRT} pour une preuve compl\`ete) \'etablit   la propagation du chaos asymptotiquement dans le scaling de Grad~: plus pr\'ecis\'ement, il montre que 
 la fonction de distribution du syst\`eme de particules est bien approch\'ee par la solution de l'\'equation de Boltzmann non lin\'eaire, avec une erreur qui tend vers 0 quand $N\to \infty$, pour un ensemble de configurations initiales dont la probabilit\'e tend vers 1 quand $N\to \infty$. 

Comme on a l'a d\'ej\`a mentionn\'e dans l'introduction, le principal d\'efaut de ce r\'esultat est qu'il n'est valable qu'en temps tr\`es court, estim\'e \`a une fraction du temps au bout duquel une particule marqu\'ee subit une collision. En particulier, il est inutilisable pour d\'ecrire les limites de relaxation rapide. Il faut noter que cette limitation n'est pas seulement li\'ee \`a la m\'ethode de preuve, mais au fait que la th\'eorie~$L^\infty$ de l'\'equation de Boltzmann est locale. D\'efinir des solutions globales pour l'\'equation de Boltzmann n\'ecessite soit d'utiliser des techniques de renormalisation \cite{diperna-lions}, soit de consid\'erer des fluctuations autour d'un \'etat d'\'equilibre.

\subsection{L'\'equation de Boltzmann lin\'eaire et le mouvement brownien}

Au niveau microscopique, la notion d'\'equilibre est reli\'ee \`a l'existence d'une mesure invariante pour le syst\`eme de $N$ sph\`eres : en notant~$Z_N := (X_N,V_N) := (x_1,\dots,x_N,v_1,\dots,v_N) \in \T_\lambda^{dN} \times \R^{dN} $, o\`u~$
	\T_\lambda := \R/ (\lambda \Z) 
	$
pour~$\lambda\geq 1$, on d\'efinit la mesure de Gibbs 
 $$M_{N,\beta} (Z_N) := { \bar \cZ_{N}^{-1} }  { \indc_{ {\mathcal D}_\eps^N } (Z_N) }\prod_{i=1}^ N
 M_\beta(v_i) \,, \quad \mbox{avec} \quad  M_\beta(v) :=  \left(\frac\beta{2\pi}\right)^{\frac{d}2 } \exp (- \frac\beta 2   |v|^2 )
 $$
 o\`u  $\indc_{ {\mathcal D}_\eps^N }$ code l'exclusion
$${\mathcal D}_\eps^N:= \big\{Z_N \in \T_\lambda^{dN} \times \R^{dN} \, / \,  \forall i \neq j \, ,\quad |x_i - x_j| > \eps \big\}$$
et o\`u~$ \bar \cZ_{N}$ est le facteur de normalisation:
$$
\bar \cZ_{N}:=  \int \indc_{ {\mathcal D}_\eps^N } (X_N) \, dX_N \, .
$$

Une premi\`ere notion de fluctuation, assez naturelle du point de vue  probabiliste, consiste \`a \'etudier la {\bf dynamique d'une particule marqu\'ee} dans le gaz globalement \`a l'\'equilibre.
Plus pr\'ecis\'ement, on consid\`ere un syst\`eme de $N$ sph\`eres dures sur  $\T_{\lambda}^d $, distribu\'ees initialement   selon la probabilit\'e
\begin{equation}
\label{donnee Boltzmann lineaire}
 f_N^0 (Z_N) :=  {   \cZ_{N}^{-1} }  { \indc_{ {\mathcal D}_\eps^N } (Z_N) } \varphi_\lambda^0(Z_1)\prod_{i=1}^ N
 M_\beta(v_i)    \, , \quad \mbox{avec} \quad \cZ_{N} := \int  { \indc_{ {\mathcal D}_\eps^N } (Z_N) } \varphi_\lambda^0(Z_1)\prod_{i=1}^ N
 M_\beta(v_i)  \, dZ_N
\end{equation}
o\`u~$\varphi_\lambda^0$ est une fonction Lipschitzienne de~$\T_{\lambda}^d \times \R^d$ telle que 
 \begin{equation}
 \label{borne-init}
  \varphi_\lambda^0(z) \leq  \mu_\lambda \, , \qquad   \int_{\T_{\lambda}^d \times \R^d}M_\beta(v) \varphi_\lambda^0(z) \, dz =1 \, ,
\end{equation}
 ce qui signifie que l'on a sp\'ecifi\'e l'\'etat de la particule 1 \`a l'instant initial.

Dans ce r\'egime, on s'attend \`a ce que la densit\'e moyenne des particules ne soit pas affect\'ee par la perturbation initiale (avec une erreur en $O(1/N)$), et que les collisions de la particule 1 soient r\'egies par un processus lin\'eaire. La distribution $M_\beta \varphi$ de la particule 1 doit alors satisfaire approximativement {\bf l'\'equation de Boltzmann lin\'eaire}
\cite{BLLS,LS2}:
\begin{equation}
\label{L-boltz}
 \begin{array}{rl}
\displaystyle  \d_t \varphi + v\cdot \nabla_x \varphi &\displaystyle =   {\mathcal L}\varphi \, , \\
  \displaystyle  {\mathcal L} \varphi (v) &:=\displaystyle 
  {1\over M_\beta(v)}  \int \Big( M_\beta(v'_1) M_\beta(v') \varphi_\lambda (v') -M_\beta(v_1) M_\beta(v)  \varphi_\lambda(v)\Big) (( v-v_1)\cdot \omega)_+ \,  d\omega dv_1  
   \end{array}
\end{equation}
L'article \cite{BGSR1} montre que cette approximation est valable sur un temps suffisamment long pour \'etudier la limite de diffusion.

\bigskip
\begin{theoreme}\label{thm1}
Soit la dynamique sur  $\T_{\lambda}^d \times \R^d$ de $N$ sph\`eres dures distribu\'ees initialement   selon (\ref{donnee Boltzmann lineaire})-(\ref{borne-init}).
  \begin{itemize}
 \item Alors, la distribution $\lambda^d f_N^{(1)}(\lambda^2 \tau , \lambda y,v)$ de la particule  marqu\'ee
 est proche de la solution de l'\'equation de Boltzmann lin\'eaire
\begin{equation}
\label{boltzlin}
  \d_\tau \varphi_\lambda +\lambda v\cdot \nabla_y \varphi_\lambda  =  \lambda^2  {\mathcal L}\varphi_\lambda  \hbox{ dans } \R^+ \times \T^d\times \R^d \end{equation}
 de donn\'ee initiale $ \varphi_\lambda^0 
$, au sens o\`u pour tout $\delta > 1$, on a dans la limite $N\eps^{d-1}\lambda ^{-d} \equiv 1$
\begin{equation}
\label{eq: approximation Boltzmann lineaire}
\big\|  \lambda^d f_N^{(1)} (\lambda^2 \tau , \lambda y,v)- M_\beta \varphi_\lambda( \tau ,   y,v) \big\|_{L^\infty([0, T]\times \T^d\times \R^d)} \leq C   \mu_\lambda  
\frac{  (\lambda^2 T)^\delta}{ ({\log\log N})}  \, .
\end{equation}
 \item En particulier, si
 $$ \lim_{\lambda \to \infty} 
  \lambda^d  \int M_\beta(v) \varphi_\lambda^0(  \lambda  y,v) \, dv
 = \rho_0(y) \hbox{ dans } L^1(\T^d)\,,$$
la distribution $\lambda^d f_N^{(1)}(\lambda^2 \tau , \lambda y,v)$ converge  au sens des mesures, quand
 $\lambda$ tend vers l'infini (beaucoup moins vite que $( \log \log N)^{1/2}$),
%~$\lambda\to \infty$,  
vers la solution $\rho(\tau, y) M_\beta(v) $ de l'\'equation de la chaleur lin\'eaire
\begin{equation}
\begin{array}l
\displaystyle \d_\tau \rho -\kappa \Delta_y \rho = 0 \hbox{ dans } \R^+ \times \T^d, \qquad \rho_{| \tau =0} (y) = \rho_0 (y)%\\
%\displaystyle  
\quad {\rm avec} \quad
\kappa: =   \int_{\R^d}  v \mathcal{L}^{-1} v \; M_\beta(v) dv \,.
\end{array}
\end{equation}
 Et le processus
$
\displaystyle \Xi (\tau) = \frac{1}{\lambda } x_1 \big( \lambda ^2 \tau \big)  
$  converge  en loi vers un mouvement brownien de variance~$\kappa$
 initialement distribu\'e sous la mesure $\rho_0$.
\end{itemize}
\end{theoreme}

\medskip
\noindent 
La limite de diffusion de l'\'equation de Boltzmann lin\'eaire~(\ref{boltzlin}) ne pose pas de difficult\'e particuli\`ere \cite{BSS,KLO}~: l'\'equation de la chaleur s'obtient  par d\'eveloppement asymptotique, et la convergence par une simple estimation d'\'energie.  La principale difficult\'e ici est donc d'obtenir l'approximation du syst\`eme de sph\`eres dures par l'\'equation de Boltzmann lin\'eaire en temps long.  A noter qu'ici, contrairement aux \'etudes de la diffusion \`a partir du gaz de Lorentz \cite{S2}, on suppose que toutes les particules ont une dynamique et interagissent entre elles. 

Le point cl\'e de la preuve consiste \`a montrer qu'on peut ``it\'erer" l'argument de Lanford, en utilisant les contr\^oles globaux fournis par le principe du maximum et une comparaison avec la mesure invariante.

\smallskip \noindent Notons que le th\'eor\`eme peut ais\'ement \^etre \'etendu au cas o\`u la donn\'ee initiale~(\ref{donnee Boltzmann lineaire}) s'\'ecrit
 $$
 f_N^0 (Z_N) :=  {   \cZ_{N}^{-1} }  { \indc_{ {\mathcal D}_\eps^N } (Z_N) } \varphi_\lambda ^0(Z_s)\prod_{i=1}^ N
 M_\beta(v_i) \, ,
$$
avec~$s$ fini, ce qui correspond \`a modifier la distribution d'un nombre fini de particules.

  \subsection{G\'en\'eralisation \`a l'\'equation de  Boltzmann lin\'earis\'ee et les \'equations de Stokes incompressibles}
 Une g\'en\'eralisation naturelle du th\'eor\`eme pr\'ec\'edent serait d'obtenir les \'equations de Stokes incompressibles en passant par   l'\'equation de  Boltzmann lin\'earis\'ee.
  Avec la notion pr\'ec\'edente de fluctuation, seul un tr\`es petit nombre de particules voient leur densit\'e affect\'ee par la perturbation initiale~: en particulier, cette perturbation n'a aucun effet sur la densit\'e moyenne. Par cons\'equent, on ne suit pas la r\'etroaction des particules marqu\'ees sur le fond, et la dynamique macroscopique se r\'eduit \`a une \'equation sur la densit\'e.

Afin d'obtenir   l'\'equation de  Boltzmann lin\'earis\'ee
\begin{equation}
\label{LL-boltz}
 \d_t \varphi +  v\cdot \nabla_x \varphi  =  \int M_\beta(v_1) ( \varphi (v')+\varphi(v'_1) -\varphi(v)-\varphi(v_1)) (( v-v_1)\cdot \omega)_+ d\omega dv_1 
\end{equation}
il faut perturber la mesure d'\'equilibre de fa\c con sym\'etrique, avec des {\bf fluctuations de tr\`es petite taille}. On se ram\`ene alors  apr\`es changement d'\'echelle \`a \'etudier  la dynamique sur  $\T_{\lambda}^d \times \R^d$ de $N$ sph\`eres dures distribu\'ees initialement   selon 
$$ f_N^0 (Z_N) :=  {   \cZ_{N}^{-1} }  { \indc_{ {\mathcal D}_\eps^N } (Z_N) } \Big(\sum_{i=1}^N \varphi^0 (z_i) \Big)  \prod_{i=1}^ N
 M_\beta(v_i)   \hbox{ de sorte que } \int{ (f_N^0)^2\over M_N} dZ_N \leq C_\lambda N\,.
$$
La limite de diffusion de l'\'equation de Boltzmann lin\'earis\'ee  vers la solution $ M_\beta (\rho +u\cdot v +\theta {|v|^2-3\over 2})$ des \'equations de Stokes incompressibles
$$\begin{array}l
\displaystyle \d_\tau\theta -\kappa \Delta_y \theta = 0 ,\quad \nabla_y (\rho+\theta) =0,\\
\displaystyle \d_\tau u -\mu \Delta_y u +\nabla_y p =0,\quad \nabla_y \cdot u =0,
\end{array}\quad  \hbox{ sur } \R^+ \times \T^d
$$
est maintenant bien comprise \cite{BGL2}. En revanche la limite du syst\`eme de sph\`eres dures vers  l'\'equation de Boltzmann lin\'earis\'ee reste un probl\`eme en cours d'\'etude~\cite{BGSR2}, en raison de la difficult\'e \`a  \'etendre l'analyse de Lanford dans un cadre fonctionnel qui est typiquement $L^2$ et non plus $L^\infty$, dans lequel la notion de trace est beaucoup moins naturelle.

  \begin{remarque}~: Il est important de noter (voir par exemple \cite{BLLS}) que l'\'equation de Boltzmann non lin\'eaire (\ref{NL-boltz}), l'\'equation de Boltzmann lin\'eaire (\ref{L-boltz}) et l'\'equation de Boltzmann lin\'earis\'ee (\ref{LL-boltz}) sont trois fermetures diff\'erentes de la m\^eme hi\'erarchie infinie, obtenue en passant formellement \`a la limite dans la hi\'erarchie BBGKY
(\ref{BBGKY}). La relation de fermeture est enti\`erement prescrite par la donn\'ee initiale.
\end{remarque}
    
\section{Recollisions et propagation du chaos}

\subsection{Repr\'esentation par des arbres de collisions}
La preuve originelle de Lanford \cite{lanford} est assez descriptive. L'id\'ee est de caract\'eriser la position et la vitesse d'une particule \`a un instant donn\'e en d\'ecrivant son arbre de collisions jusqu'\`a cet instant $t$, et la configuration initiale de toutes les particules appartenant \`a cet arbre.

Le point de d\'epart est l'\'equation de Liouville qui traduit les \'equations de la dynamique microscopique en terme de probabilit\'e sur l'espace des configurations \`a $N$ particules. Pour des sph\`eres dures, elle s'\'ecrit simplement
$$ \d_t f_N +\sum_i v_i\cdot \nabla_{x_i} f_N=0 \hbox{ sur } 
{\mathcal D}_\eps^N
$$
avec  condition de r\'eflexion sp\'eculaire sur le bord du domaine (correspondant \`a un rebond \'elastique).
La formule de Green montre alors que pour tout $ s<N$, la marginale d'ordre $s$ 
$$
 f_N^{(s)} (t,Z_s ):=\int  f_N  (t,Z_N )  { \indc_{ {\mathcal D}_\eps^N } (Z_N) } dz_{s+1} \dots dz_N
$$satisfait l'\'equation
\begin{equation}
\label{BBGKY}
	(\d_t +\sum_{i=1}^s v_i\cdot \nabla_{x_i} ) f_N^{(s)} (t,Z_s) = \big( C_{s,s+1} f_N^{(s+1)}\big) (t,Z_s)\hbox{ sur } \cD_\eps^s, 
	\end{equation}
		o\`u le terme de collision est d\'efini par
		$$
	\big( C_{s,s+1} f_N^{(s+1)}\big) (Z_s) := ( N-s) \eps^{d-1} \sum_{i=1}^s \int_{{\mathbf S}^{d-1} \times \R^ d}  f_N^{(s+1)}(Z_s, x_i+\eps \omega, v_{s+1})  (v_{s+1}-v_i) \cdot \omega  \, d\omega dv_{s+1}\,.	$$
	 Ce syst\`eme de $N$ \'equations coupl\'ees est appel\'e hi\'erarchie BBGKY (d'apr\`es Bogoliubov, Born, Green, Kirkwood, Yvon  \cite{lanford,CIP}). 
% \cite{bogoliubov,born,Kirkwood,Yvon}).
La premi\`ere  \'equation ($s=1$) est assez semblable \`a l'\'equation de Boltzmann, mais elle ne se d\'ecouple du reste de la hi\'erarchie que dans le cas o\`u {\bf les particules collisionnelles sont ind\'ependantes}.

La repr\'esentation par des arbres de collision s'obtient en introduisant les op\'erateurs de transport~$ {\bf S}_s$ \`a~$s$ particules dans~$\cD_\eps^s$ et de scattering, et en int\'egrant la  hi\'erarchie avec la formule de Duhamel 
$$  f^{(s)} _N(t) =  \sum_{n=0}^{N-s} Q_{s,s+n} (t) f^{(s+n)}_N(0)$$
o\`u l'on a d\'efini 
$$
Q_{s,s+n} (t) :=   \int_0^t \int_0^{t_1}\dots  \int_0^{t_{n-1}}  {\bf S}_s(t-t_1) C _{s,s+1}  {\bf S}_{s+1}(t_1-t_2) 
  C _{s+1,s+2}\dots   {\bf S}_{s+n}(t_n)   \: dt_{n} \dots dt_1\,.
	$$
Chaque terme de la somme peut alors  \^etre associ\'e \`a une pseudo-dynamique d\'efinie de fa\c con r\'ecursive
\begin{itemize}
\item les $s+i$ particules sont transport\'ees (dans le sens des temps d\'ecroissants) entre  $t_i$ et $t_{i+1}$;
\item une particule suppl\'ementaire est  ``adjointe" au syst\`eme au temps $t_{i+1}$.
\end{itemize}

\begin{figure}[h] %  figure placement: here, top, bottom, or page
   \centering
  \includegraphics[width=3in]{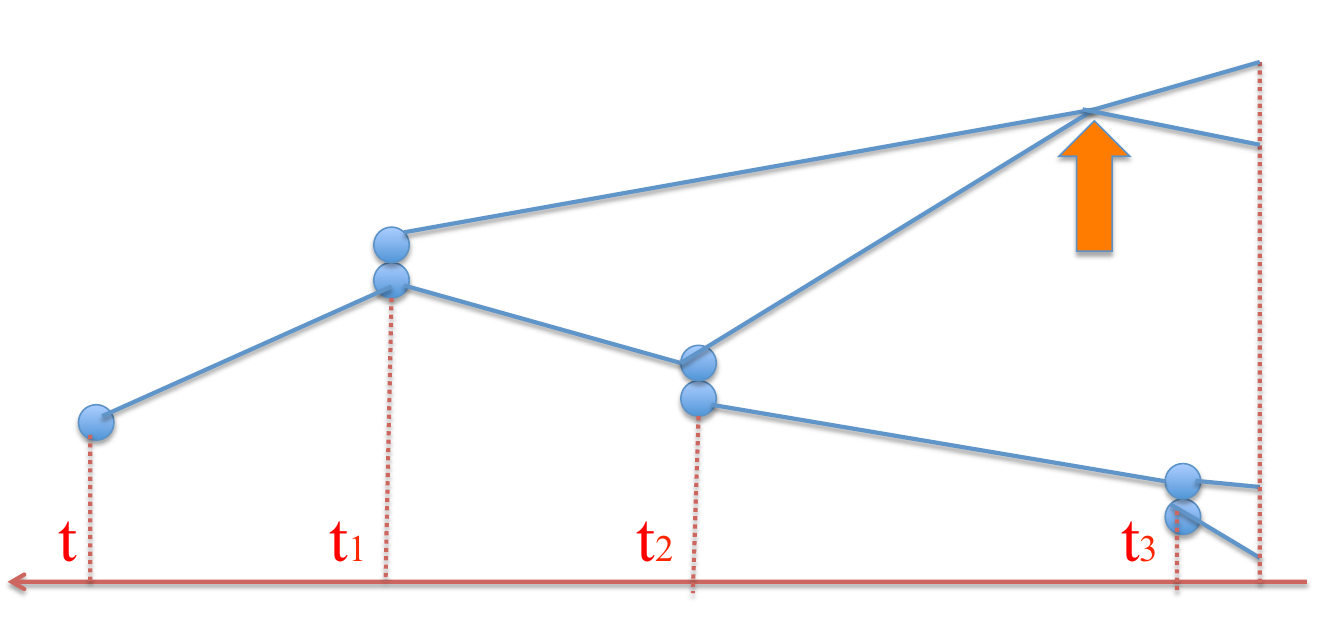} \includegraphics[width=3in]{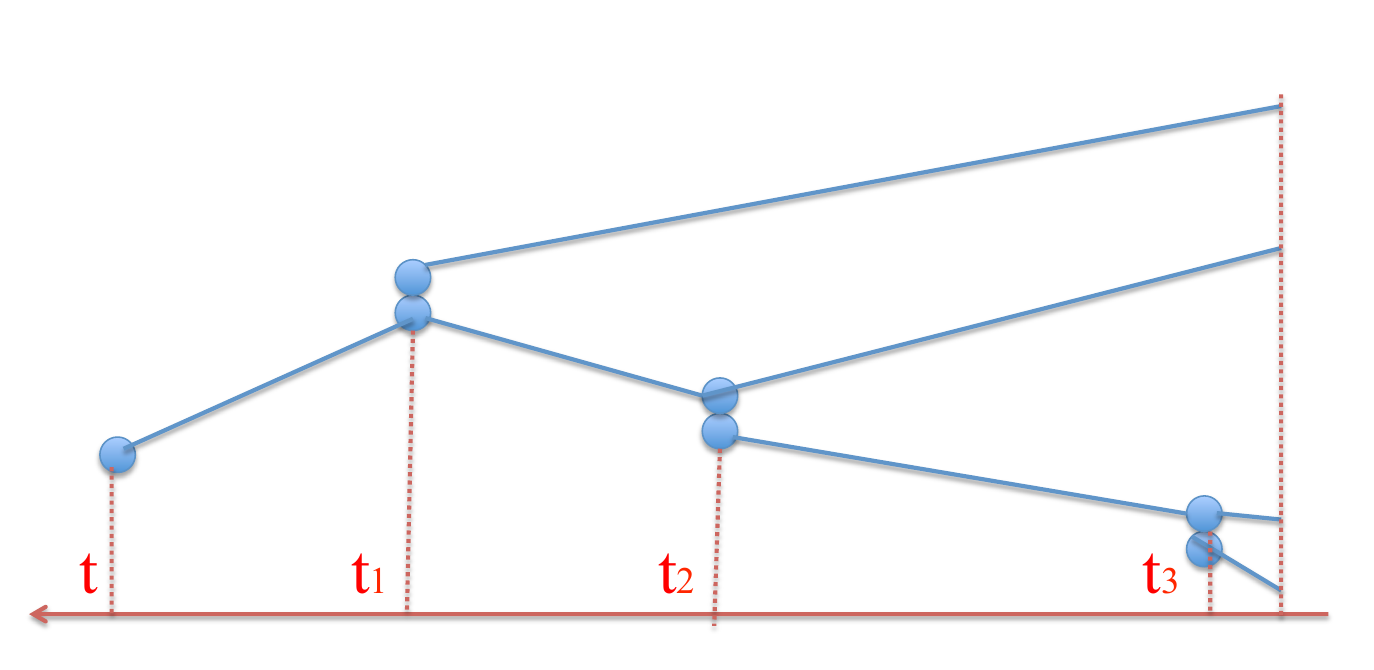} 
  \caption{Arbre de collision avec et sans recollision }
 \end{figure}

Les op\'erateurs de transport ${\bf S}_{s+i}$ co{\"\i}ncident avec le transport libre si et seulement si il n'y a pas de recollision.
Autrement dit, les solutions de l'\'equation de Boltzmann donnent une bonne approximation de la hi\'erarchie BBGKY  tant qu'il n'y a pas recollision.
La preuve de convergence consiste donc \`a montrer que les recollisions sont de probabilit\'e asymptotiquement nulle quand $\eps \to 0$.

\subsection{Lemmes g\'eom\'etriques}

La d\'emonstration propos\'ee dans \cite{GSRT} est it\'erative~:  on d\'efinit  les ``bonnes configurations"
$$ \left\{ Z_k \in \T_\lambda^{dk} \times \R^{dk}\,/\, \forall s \in [0,t]  ,\quad \forall i\neq j ,
\quad  d(x_i -s v_i, x_j- sv_j) \geq   \eps_0\right\}\,,$$
 o\`u $d$ est  la distance sur le tore $\T_\lambda^d$ et $\eps_0 \gg \eps$. Et on montre que ces bonnes configurations
sont stables par adjonction d'une particule pourvu qu'on exclue un petit ensemble de vitesses et d'angles de d\'eflection dits pathologiques.

Une \'etude g\'eom\'etrique du transport libre permet de contr\^oler la taille de cet ensemble en fonction de~$\eps$ et de $s$. L'estimation \'el\'ementaire pour chaque paire de particules est la suivante~:

\smallskip

\begin{lemme}
 Soient deux positions~$  x_1,  x_2 \in \T_\lambda^d$ telles que $d(   x_1,  x_2)\geq \eps_0 \gg \eps$,  et une vitesse~$ v_1$ telle que~$|v_1| \leq E < \infty$. Pour~$\delta,t>0$ donn\'es, il existe un ensemble~$K(  x_1-  x_2) $ de mesure petite
 $$ 
 |K(  x_1-  x_2) | \leq  CE^d \left(  \left({ \eps \over \eps_0 }\right)^{d-1}+  \left({ \eps_0 \over E\delta}\right)^{d-1}+ \Big(\frac{E t}\lambda \Big)^d \Big( \frac{\eps_0 }\lambda\Big) ^{d-1}\right) $$
  tel que pour toute vitesse~$  v_2 \notin (v_1+ K(  x_1-  x_2))$, avec~$|v_2| \leq E $,
 \begin{itemize}
 \item[(i)]  il n' y a {\bf pas de collision  sur  $[0,t]$} par le flot inverse:  pour tout~$s>0$,
  $d(x_1-v_1 s, x_2-v_2 s) >\eps$;
     \item[(ii)]  on est \`a nouveau dans une {\bf bonne configuration apr\`es un temps $\delta$}: pour tout~$s>\delta$, 
   $d(x_1-v_1 s, x_2-v_2 s)>\eps_0 .$
   \end{itemize}

\end{lemme}

Dans le cas de l'espace entier, l'ensemble $K$ est essentiellement la r\'eunion d'un c\^one et d'un cylindre (intersect\'es avec la boule de rayon $E$).
Ici \`a cause de la p\'eriodicit\'e du domaine, il appara\^\i t un troisi\`eme terme dans l'estimation.
On notera que ce lemme n'a pas d'analogue simple dans le cas des domaines \`a bords.

\begin{figure}[h] %  figure placement: here, top, bottom, or page
   \centering
  \includegraphics[width=3.5in]{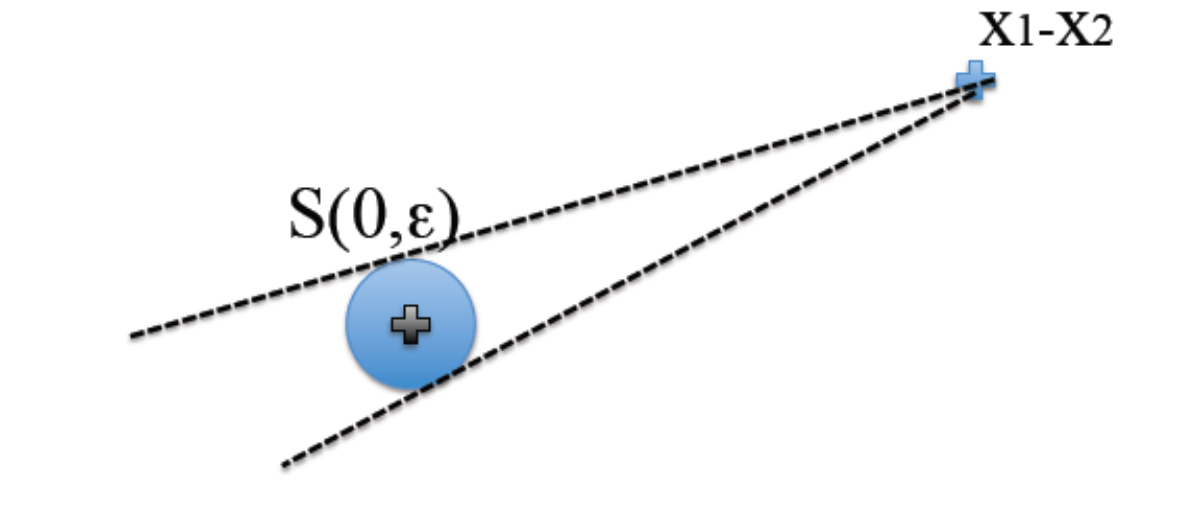} 
  \caption{C\^one des vitesses donnant une recollision (violation de la condition (i))}
 \end{figure}

Pour redresser les trajectoires et se ramener au transport libre, on modifie alors successivement   chaque op\'erateur de collision en enlevant l'ensemble pathologique du domaine d'int\'egration~: on doit donc r\'ep\'eter l'op\'eration environ $n(s+n)$ fois pour l'op\'erateur~$Q_{s,s+n}$.

\section{Contr\^ole du processus de branchement pour une particule marqu\'ee}

Cet algorithme de redressement n'est admissible que si les  arbres  de collision sont assez petits (avec un nombre $n$ de points de branchement au plus de l'ordre de  $\log N$).
Pour des donn\'ees initiales assez g\'en\'erales, la taille des arbres est contr\^ol\'ee  par une estimation de type Cauchy-Kowalewski pour la hi\'erarchie BBGKY, ce qui explique qu'on ne peut atteindre que des temps tr\`es courts. Cette difficult\'e est contourn\'ee dans \cite{BGSR1} en consid\'erant des donn\'ees initiales proches de l'\'equilibre, au sens de (\ref{donnee Boltzmann lineaire}).

\subsection{Mesure invariante et principe du maximum}

 Sous l'hypoth\`ese (\ref{donnee Boltzmann lineaire}), on peut en effet obtenir des bornes a priori globales en temps et essentiellement uniformes en $N$ en utilisant la comparaison avec la mesure invariante 
 $$ f_N^0(Z_N) \leq \mu_\lambda    M_{N,\beta}$$
 et le principe du maximum. On en d\'eduit que 
$$ \sup_{t }  f_N^{(s)}(t,Z_s) \leq    \mu_\lambda   M_{N,\beta}^{(s)} (Z_s)\leq  \lambda^{-ds}     \mu_\lambda    \big( 1 - C\eps  \big)^{-s}M_\beta^{\otimes s} (V_s)\,.$$
Ces bornes permettent de contr\^oler la taille des arbres, et plus pr\'ecis\'ement de montrer que les arbres de grande taille sont de probabilit\'e asymptotiquement nulle.

\subsection{Arbres \`a croissance super-exponentielle}

On dit qu'un arbre de  collision est admissible s'il a moins de $n_k= a^{ k}$  points de branchement sur l'intervalle~$[t-k\tau ,t-(k-1) \tau]$, o\`u~$\tau$ est un param\`etre \`a choisir suffisamment petit et $a$ est une constante \`a choisir suffisamment grande.

\begin{figure}[h] %  figure placement: here, top, bottom, or page
   \centering
  \includegraphics[width=3.5in]{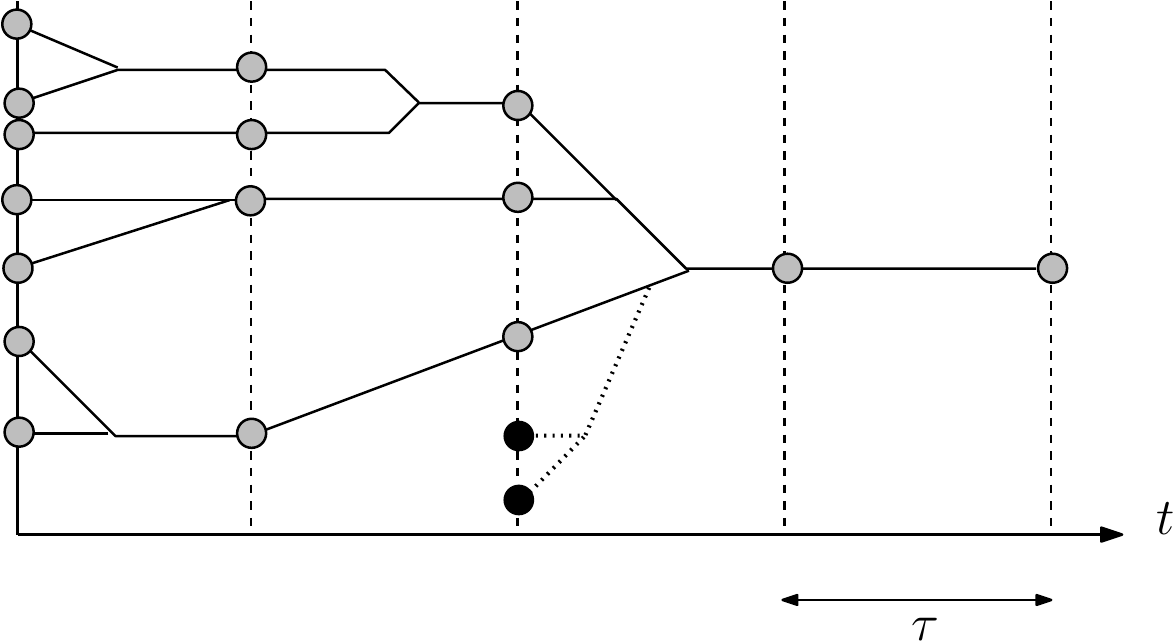} 
  \caption{Arbres de collision admissible (gris) et non admissible (avec l'ajout des particules noires)}
 \end{figure}

On introduit par cons\'equent la d\'ecomposition
$$\begin{array}l
\displaystyle  f^{(1,K)}_N(t)    : =   \sum_{j_1=1}^{n_1-1}\! \!   \dots \!  \! \sum_{j_K=0}^{n_K-1}Q_{1,J_1} (\tau )Q_{J_1,J_2} (\tau )
 \dots  Q_{J_{K-1},J_K} (\tau ) \, f^{0(J_K)}_N  \\
  \displaystyle R_N^K(t)  := \sum_{k=1}^K \; \sum_{j_1=1}^{n_1-1} \! \! \dots \! \! \sum_{j_{k-1}=0}^{n_{k-1}-1} \; 
Q_{1,J_1} (\tau ) \dots  Q_{J_{k-2},J_{k-1}} (\tau ) \, R_{J_{k-1},n_k}(t-k \tau, t-(k-1)\tau )  \, . 
\end{array}$$
o\`u on note $J_k = \sum_{i=1} ^k j_i$.

La brique \'el\'ementaire ici est l'estimation de continuit\'e \`a perte pour les op\'erateurs de collision
$$\| Q_{s,s+n}(t)f_{s+n} \|_{\eps,s,\frac\alpha2} \leq    e^{s-1} \left( {C_d\lambda ^d t\over \alpha^{\frac {d+1}2}}  \right) ^n 	\| f_{s+n}	\| _{\eps,s+n,\alpha} $$
dans les espaces \`a poids $X_{\eps, s, \alpha}$ d\'efinis par
$$ 	\| f_{s}	\| _{\eps,s,\alpha} := \sup _{Z_s \in \cD_\eps^s} \Big(\exp (\frac\alpha 2 |V_s|^2) |f_s(Z_s)| \Big) \,.$$
En choisissant $\tau =\gamma/t$, 
on obtient  le contr\^ole suivant sur le terme de reste
$$ \lambda^d\left\|
  R_N^K(t)\right\|_{L^\infty({\mathbf T}_\lambda^{d} \times \R^{d})}   \leq  C    \mu_\lambda  \gamma^a   \, .$$
En appliquant cette m\^eme estimation de continuit\'e \`a perte pour contr\^oler les erreurs it\'er\'ees dues au redressement des trajectoires, et en optimisant le param\`etre $\gamma$, on obtient l'in\'egalit\'e (\ref{eq: approximation Boltzmann lineaire}) 
 pour des temps minor\'es par $(\log \log N)^\frac{a -1}{a}$ o\`u $a$ peut \^etre choisi arbitrairement grand.

\smallskip
\begin{remarque}~:
Notons que la convergence dans le th\'eor\`eme \ref{thm1} est \'etablie pour des temps  presque de l'ordre de $\log \log N$. Cette limitation sur le temps semble optimale
pour cette m\'ethode, dans la mesure o\`u  les arbres de collision croissent exponentiellement en~$t$, et que  $\log N$ est la taille critique des arbres au-del\`a de laquelle des transitions de phase peuvent apparaitre (on s'attend alors \`a ce que des recollisions aient lieu avec probabilit\'e strictement positive).
\end{remarque}

\begin{remarque}~:
L'approche du mouvement brownien que nous avons pr\'esent\'ee  s'inscrit dans la lign\'ee des travaux en physique qui ont permis de mettre en \'evidence le mouvement Brownien \`a partir d'un syst\`eme m\'ecanique \cite{duplantier,kahane}. 
%\cite{smoluchowski,langevin,duplantier,kahane}. 
Dans ces travaux, la particule marqu\'ee (particule de pollen) est de masse tr\`es sup\'erieure aux particules l'environnant. Son mouvement est donc r\'egi par un grand nombre de d\'eflections infinit\'esimales.

La situation que l'on consid\`ere ici est un peu diff\'erente puisque la particule marqu\'ee est de m\^eme taille que les autres et subit une d\'eflection macroscopique \`a chaque impact.
N\'eanmoins les m\'ethodes utilis\'ees dans \cite{BGSR1} pour l'analyse asymptotique semblent assez robustes pour pouvoir \^etre \'etendues \`a des
dynamiques o\`u la particule marqu\'ee aurait une taille et une masse diff\'erente de celles des autres particules.
%la dynamique du grain de pollen.
\end{remarque}

\end{document}